# On estimating the change point in generalized linear models[*]

## Hongling Zhou[1] and Kung-Yee Liang[2]

*US Food and Drug Administration and Johns Hopkins University*

**Abstract:** Statistical models incorporating change points are common in practice, especially in the area of biomedicine. This approach is appealing in that a specific parameter is introduced to account for the abrupt change in the response variable relating to a particular independent variable of interest. The statistical challenge one encounters is that the likelihood function is not differentiable with respect to this change point parameter. Consequently, the conventional asymptotic properties for the maximum likelihood estimators fail to hold in this situation. In this paper, we propose an estimating procedure for estimating the change point along with other regression coefficients under the generalized linear model framework. We show that the proposed estimators enjoy the conventional asymptotic properties including consistency and normality. Simulation work we conducted suggests that it performs well for the situations considered. We applied the proposed method to a case-control study aimed to examine the relationship between the risk of myocardial infarction and alcohol intake.

## 1. Introduction

The problems of detecting abrupt changes at unknown points and estimating the locations of changes are known as the change point problem. The change point problem occurs frequently in medical research. For example, cancer incidence rates remain relatively stable for people at a younger age, but change drastically after a certain age threshold (MacNeill and Mao [13]). The data obtained from a group of preschool boys indicates that their weight/height ratio relates to their age in one way before a certain age but that the functional relation of the two changes afterwards (Gallant [8]). Another example arises from a study of the risk of myocardial infarction(MI), which showed a sharp decrease in risk at low alcohol intakes and a dramatic increase after reaching a certain amount of daily alcohol consumption (Pastor and Guallar [16]). Although these examples each have distinctive features of their own, the common theme here is that the relationship of the response variables and a covariate of interest is subject to an abrupt change at a certain threshold. Very often, scientists are interested in the threshold for clinical or preventive purposes. The change point model is useful in that it purposely includes a parameter to capture the notion of threshold.

[*]Supported in part by NIH Grant GM49909.

[1]10903 New Hampshire Ave, Silver Spring, MD, 20903, USA, e-mail: hongling.zhou@fda.hhs.gov
[2]615 North Wolfe Street, Baltimore, MD 21205, USA, e-mail: kyliang@jhsph.edu

In this paper, we focus on the estimation of the change point in the so-called broken-line regression models, where the regression function is assumed to be continuous at the point of change. Estimation for these type of models with normally distributed responses has been developed by various authors. Sprent [18] was among the first to discuss the estimation of the piecewise linear models. His interest in this type of model is based on the observation that a biologist would often postulate a two-phase linear model over some alternatives such as the quadratic model largely on intuitive grounds. Hinkley [10] considered the same two-phase straight line model and derived the maximum likelihood estimator (MLE) of the change point by its marginal likelihood function and presented the asymptotic distribution of the estimator. Feder [7] studied the model in a more general framework and proved the consistency of the least-squares estimators of the regression coefficients and the change point. The estimators are asymptotically normal for some special cases including models with all linear segments. Bhattacharya [2, 3] presented the asymptotic properties of the change point and regression coefficient estimators using a local log likelihood process approach. Through this approach, he showed the distinctive features of the asymptotic properties of the change point with and without the continuity constraint at the point of change.

A major difficulty in estimating the change point as a parameter for regression models is the nonsmoothness of the likelihood function with respect to the change point considered as a parameter. Many authors tried to circumvent this problem by using various smooth transitions between the two linear regimes, or using other types of the function such as the quadratic function in one of the segments separated by the change point. This technique was mostly used for the model with normally distributed response variables. Gallant and Fuller [9] discussed such a model and used a modified Gauss-Newton method to obtain the least squared estimates. The simulation studies showed that the estimates obtained by this method work well, but no asymptotic properties were presented in their work. Bacon and Watts [1] proposed a model which can accommodate a smooth transition as well as an abrupt change with a Bayesian estimation procedure to determine the parameter values. To the best of our knowledge, with the exception of normal distributions, there is very little research being done where the response variables are categorical. Furthermore, the asymptotic properties of the estimators are virtually unavailable.

In Section 2, we propose an estimation procedure for the change point and the corresponding regression coefficients in the framework of generalized linear models, and present the asymptotic properties of the proposed estimators. We discuss in Section 3 some extensions of the estimation procedure to more general models involving a change point. Section 4 discusses the choice of the bandwidth parameter and the smoothing function for the proposed estimation procedure. Section 5 presents the simulation results to assess the finite sample performance of the proposed method. In Section 6, we apply the proposed method to a data example from the EURopean study on Antioxidant, Myocardial Infarction, and breast Cancer (EURAMIC, Kardinaal et al. [12]) followed by discussion in Section 7. The proofs of the propositions are given in the Appendices.

## 2. Proposed estimation procedure for the change point

Our approach is closely related to the density estimation, whereas the cumulative distribution function (CDF) is estimated by the empirical distribution function (EDF), but the density function cannot be estimated by the derivative of CDF, which in essence is a step function (Rao [17]). In density estimation, a kernel function (usually a known density function itself) is introduced to solve this problem by making the EDF smooth and differentiable everywhere. The proposed estimation procedure is directly motivated by the approach adopted by Horowitz [11] and in the Econometrics literature. In particular, Manski [14] showed that the maximum score estimator of the coefficient vector of a binary response model is consistent under a weak distributional assumption. However, the asymptotic distribution of the maximum score estimator was not presented since the objective function it maximizes is a step function. Horowitz [11] used a smooth version of the same objective function to make it continuous and differentiable. The procedure produced an estimator that not only converges more rapidly but also has a tractable asymptotic distribution.

For a sequence of random variables $Y_i, i = 1, \ldots, n$, having probability (density) function of the form

$$(2.1) \qquad f(y_i; \theta) = c(y_i) exp(\theta y_i - b(\theta)), \quad c \geq 0$$

with the natural link function, we consider the following model incorporating the unknown change point $\tau$ as

$$(2.2) \qquad \theta = g(\mu) = \beta_0 + \beta_1 x + \beta_2 (x - \tau)_+,$$

where $a_+ = aI(a > 0)$ and $I(a > 0)$ is the indicator function, and $\beta_2 \neq 0$ for identifiability of $\tau$ (Davies [5]). Here $\beta_0$ and $\beta_1$ are the intercept and the slope relating the response variable $Y$, through the link function $g$, to the covariate $x$ for $x \leq \tau$, and $\beta_2$ is the difference in slopes for the segments before and after the change point $\tau$. With the traditional likelihood approach, the likelihood function is not differentiable with respect to the change point $\tau$. Specifically, the indicator function $I(x > \tau)$ is not differentiable with respect to $\tau$. Consequently, one of the regularity conditions for the usual asymptotic theory, namely, a certain degree of smoothness of the objective function with respect to the parameters, is violated. To circumvent this critical problem, define a continuous function $K(\cdot)$ which satisfies

1. $|K(u)| < M$ for some M, $0 < M < \infty$ and $u \in (-\infty, \infty)$.
2. $\lim_{u \to -\infty} K(u) = 0$ and $\lim_{u \to \infty} K(u) = 1$.

We propose to estimate the change point $\tau$ as well as the regression coefficients $\beta_0, \beta_1$ and $\beta_2$ by maximizing the following objective function:

$$\begin{aligned} Q_n(\beta, \tau) &= \sum_{i=1}^{n} [y_i \{ \beta_0 + \beta_1 x_i + \beta_2 (x_i - \tau) K(\frac{x_i - \tau}{h_n}) \} \\ &\qquad - b(\beta_0 + \beta_1 x_i + \beta_2 (x_i - \tau) K(\frac{x_i - \tau}{h_n}))] \\ &= \sum_{i=1}^{n} q_i(\beta, \tau), \end{aligned}$$

where $\{ h_n : n = 1, 2, \ldots \}$ is a sequence of positive numbers satisfying $\lim_{n \to \infty} h_n = 0$. Here $K(\cdot)$ is analogous to a cumulative distribution function rather than a density function, the latter being more commonly used in problems such as density estimation.

Note that the difference between this objective function $Q_n$ and the otherwise usual log-likelihood function is that the indicator function $I(x > \tau)$ in the likelihood function is replaced by the smoothing function $K((x - \tau)/h_n)$. It is clear that the objective function $Q_n$ is twice differentiable with respect to all parameters, and with some suitable conditions on the distribution for $x$, $K((x - \tau)/h_n)$ converges to $I(x > \tau)$ uniformly as $n \to \infty$. As shown below, these two factors play key roles in the asymptotic behavior of the estimates which maximize $Q_n$.

Define

$$
\begin{aligned}
\delta &= (\beta, \tau) = (\beta_0, \beta_1, \beta_2, \tau), \\
S_n(\delta) &= \bigtriangledown Q_n(\delta) = \left(\frac{\partial Q_n(\beta, \tau)}{\partial \beta_0}, \frac{\partial Q_n(\beta, \tau)}{\partial \beta_1}, \frac{\partial Q_n(\beta, \tau)}{\partial \beta_2}, \frac{\partial Q_n(\beta, \tau)}{\partial \tau}\right)^t, \\
\Sigma_n(\beta, \tau) &= \mathbf{cov}_{\beta, \tau} S_n(\beta, \tau)
\end{aligned}
$$

and

$$
J_n(\beta, \tau) = -\bigtriangledown S_n(\delta).
$$

Note that we use the conditions imposed on the limiting configuration of the covariate X, similar to those of Bhattacharya [3]. Suppose that there is a function $F$ with $\alpha = F(\tau^0)$, for $x_1 \leq x_2 \cdots, \leq x_{n\alpha} \leq \tau \leq x_{n\alpha+1} \leq \cdots \leq x_n$, and as $n \to \infty$,

$$
(2.3) \qquad (n\alpha)^{-1} \sum_{i=1}^{n\alpha} \begin{pmatrix} 1 & x_i \\ x_i & x_i^2 \end{pmatrix} \to \begin{pmatrix} 1 & \mu_1 \\ \mu_1 & \mu_1^2 + \sigma_1^2 \end{pmatrix},
$$

$$
(2.4) \qquad (n - n\alpha)^{-1} \sum_{i=n\alpha+1}^{n} \begin{pmatrix} 1 & x_i \\ x_i & x_i^2 \end{pmatrix} \to \begin{pmatrix} 1 & \mu_2 \\ \mu_2 & \mu_2^2 + \sigma_2^2 \end{pmatrix}.
$$

Following Fahrmeir and Kaufmann [6], to ensure the asymptotic normality of $\hat{\delta}$, some additional assumptions are needed. Specifically, define a neighborhood of true parameter values $\delta^0 = (\beta_0^0, \beta_1^0, \beta_2^0, \tau^0)$ as

$$
N_n(\omega) = \{\delta : ||(\Sigma_n^{1/2})^t (\delta - \delta^0)|| \leq \omega\}, \quad n = 1, 2, \ldots,
$$

where $\omega > 0$, $\Sigma_n$ is simplified for $\Sigma_n(\delta^0)$, and $\Sigma_n^{1/2}$ is defined as the square root of positive definite matrix $\Sigma_n$ such that $\Sigma_n = \Sigma_n^{1/2}(\Sigma_n^{1/2})^t$. The conditions are

(I) $\lim_{n \to \infty} \frac{1}{n} \Sigma_n = \Sigma$, where $\Sigma$ is finite and positive definite.

(II) For all $\omega > 0$, $\max_{\delta \in N_n(\omega)} ||V_n(\delta) - I|| \to 0$ in probability under both measures $P_{\delta^0}$ and $P_{\delta_n}$ with $\delta_{\mathbf{n}} = \delta^0 + \omega(\Sigma_n^{-1/2})^t \lambda$, where $V_n(\delta) = \Sigma_n^{-1/2} J_n(\delta)(\Sigma_n^{-1/2})^t$ and $\lambda'\lambda = 1$.

(III) $g(\cdot)$ is twice continuously differentiable $g(\cdot), g'(\cdot)$ and $g''(\cdot)$ are bounded.

**Proposition 1** (Asymptotic normality). *Let $\hat{\boldsymbol{\delta}} = (\hat{\beta}, \hat{\tau})$ be the estimators which maximize the objective function $Q_n(\boldsymbol{\beta}, \tau)$ and $\boldsymbol{\delta^0}$ be the true value of $\boldsymbol{\delta}$, then under some regularity conditions and given that*

(i) *$X_i$ is assumed to be bounded, where $i = 1, \ldots, n$;*

(ii) *$\lim_{n \to \infty} P(|X_n - \tau| < \epsilon) = 0$ for some $\epsilon > 0$; and*

(iii) *$\lim_{n \to \infty} \frac{1}{n} \sum_{i=1}^{n} E\{log f(Y_i; \boldsymbol{\beta}, \tau, x_i); \boldsymbol{\beta^0}, \boldsymbol{\tau^0}\} < \infty$*

*the normed estimator for δ is asymptotically normal, i.e.,*

$$\Sigma_n^{-1/2} J_n (\hat{\beta}_0 - \beta_0^0, \hat{\beta}_1 - \beta_1^0, \hat{\beta}_2 - \beta_2^0, \hat{\tau} - \tau^0)^t \overset{D}{\sim} N(\mathbf{0}, \mathbf{I_4}).$$

*Or equivalently, by the results of Lemma B.1 in Appendix B.*

$$\sqrt{n}(\hat{\beta}_0 - \beta_0^0, \hat{\beta}_1 - \beta_1^0, \hat{\beta}_2 - \beta_2^0, \hat{\tau} - \tau^0)^t \overset{D}{\sim} N(\mathbf{0}, \Sigma^{-1}).$$

Proof of Proposition 1 is given in the Appendices.

## 3. Some extensions

The model (2.2) discussed so far involves two straight lines, which is characterized in the literature as broken line regression or joint point models. Straight line describes the abrupt change mechanism more distinctively than other types of models. As discussed in the introductory section, some authors (e.g., Pastor and Guallar [16]; Gallant [8]) used a quadratic-linear or linear-quadratic model to characterize the change point. For example, a linear-quadratic model is expressed as

$$(3.1) \qquad g(\mu_i) = \beta_0 + \beta_1 x_i + \beta_2 (x_i - \tau)^2 I(x_i \geq \tau).$$

A quadratic-linear model can be expressed similarly by changing the indicator function to $I(x_i < \tau)$. The advantage of this type of model, specifically, is that for estimation purposes the likelihood function has the first derivative with respect to the change point. Thus the Fisher information can be derived as the covariance of the score function. However, the usual asymptotic properties still cannot go through.

Since both linear-linear and linear-quadratic (or quadratic-linear) models face the same non-differentiability problem in $\tau$, the approach proposed in section 2 can be easily adopted and extended to the latter situation.

**Corollary 1.** *Replacing the term $x_i - \tau$ by $(x_i - \tau)^2$ in $Q_n$, and with the regularity conditions tailored to such a replacement, the resulting estimators of $(\beta_0, \beta_1, \beta_2, \tau)$ for model (3.1) are consistent and asymptotically normal.*

The proof is similar to that of the linear-linear model and is omitted. In model (2.2), we have assumed that there is only one independent variable, which involves the change point. This method can be easily extended to situations where adjustment for additional independent variables is needed. Specifically, suppose there are k additional variables $z_1, \ldots, z_k$, the systematic component part of a generalized linear model can be modified as

$$g(\mu_i) = \beta_0 + \beta_1 x_i + \beta_2 (x_i - \tau)_+ + \gamma_1 z_1 + \cdots + \gamma_k z_k.$$

The asymptotic properties of the estimators for $\boldsymbol{\gamma}$ and $\boldsymbol{\delta}$ are similar to those in model (2.2) without additional independent variables and are thus omitted.

## 4. Computational issues

The estimation of the parameters $(\beta_0, \beta_1, \beta_2, \tau)$ involves the value of the bandwidth parameter $h_n$. For a smoothed maximum estimator of a binary response model,

Horowitz [11] suggested choosing bandwidth $h_n \propto n^{-1/(2m+1)}$, where m is the order of the kernel $K'(\cdot)$ defined as the integer satisfying the following

$$(4.1) \qquad \int_{-\infty}^{\infty} v^i K'(v) dv = \begin{cases} 0, & \text{if} \quad i < m, \\ d(nonzero), & \text{if} \quad i = m. \end{cases}$$

In our asymptotic analysis, with $K(\cdot)$ appearing as a single term in some of the elements of the hessian matrix, the form of $h_n$ does not appear to affect the asymptotic properties of the estimators as long as $h_n \to 0$ when $n \to \infty$ and meets the restriction as follows.

A sufficient condition for the uniform convergence of $K((x_i - \tau)/h_n) \to I(x_i > \tau)$, is $(x_i - \tau)/h_n \to \pm\infty$ for all $i = 1, \ldots, n$. Generally, in order to satisfy the above conditions, we need to have $h_n << \min(|x_i - x_j|, 1 \leq i, j \leq n)$. Specifically, in the iteration procedure to estimate the parameters, $h_n$ needs to be small enough so that the negative hessian matrices are positive definite.

It is recognized in the density estimation literature that the choice of the smoothing function does not affect the asymptotic properties of the density estimator. In the simulation that follows, we are listing the results where $K(\cdot)$ is chosen as the cumulative distribution function of the standard normal distribution. We have also conducted the simulation with $K(\cdot)$ as the cumulative distribution function of the exponential distribution. Results concerning the distribution of the estimators and the actual estimated values are very similar for the two different $K(\cdot)$. However, for the same n, $h_n$ can be chosen as a function of n that approaches 0 in a slower rate with $K(\cdot)$ being the exponential cumulative distribution function for the algorithm to converge. This is consistent with the conditions (a) and (b) in Lemma B.1, i.e., $\{(x-\tau)/h_n\}K'((x-\tau)/h_n) \to 0$ and $\{(x-\tau)/h_n\}^2 K''((x-\tau)/h_n) \to 0$ as $n \to \infty$. Note that $x - \tau$ is bounded. This is equivalent to requiring $h_n$ to converge to 0 in a relatively slower rate than $K'(\infty)$ and $K''(\infty)$. Since the exponential density at the tail converges to 0 at a slower rate than the normal density, the corresponding $h_n$ can be chosen this way as well.

We have found in our simulations that the objective function $Q_n$ when evaluated at the true values $\beta_0^0, \beta_1^0, \beta_2^0$, a function of $\tau$, is rather flat around the true change point value $\tau^0$. In dealing with this issue, aside from having a small tolerance level ($10^{-5}$ at most) for convergence, another crucial issue in actually carrying out the estimation algorithm is to choose the initial values of the parameters $\beta_0, \beta_1, \beta_2$ and $\tau$. In our simulations in section 5 and the example in section 6, for the data generated, we first made some graphs of nonparametric models such as the LOESS model, then we chose a number of $\tau_i, i = 1, \ldots, I$, which on the graph are near the potential change point. Next we consider the GLM models. With the fixed change point valued at $\tau_i$, the corresponding $\boldsymbol{\beta}(\tau_i)$ values were then obtained using ordinary GLM fitting software. $\tau_p = \arg\max_{\tau_i} Q_n(\boldsymbol{\beta}(\tau_i), \tau_i)$ is chosen as the initial estimate for $\tau$, and $\boldsymbol{\beta}(\tau_p)$ as the initial values for $\boldsymbol{\beta} = (\beta_0, \beta_1, \beta_2)$.

After replacing the indicator function in the likelihood function with the smoothing function $K(\cdot)$, the model can be categorized as having nonlinear parameters. For easy calculation, according to section 11.4 in McCullagh and Nelder [15], model (2.2) can be approximated by

$$g(\mu_i) = \beta_0 + \beta_1 x_i + \beta_2 u_i + c v_i,$$

where $u = (x - \tau_0)K(\frac{x-\tau_0}{h_n})$, $v = -\{K(\frac{x-\tau_0}{h_n}) + \frac{x-\tau_0}{h_n}K'(\frac{x-\tau_0}{h_n})\}$ and $\tau_0$ is an initial value of $\tau$, which supposedly is close to true $\tau$.

After this approximation, the coefficient estimates as well as their covariance matrix can be easily obtained by regular software. The estimated $\tau$ is then $\hat{\tau} = \tau_0 - \hat{c}/\hat{\beta}_2$, and the approximate standard error by the delta method is $\{(1/\hat{\beta}_2, -\hat{c}/\hat{\beta}_2^2)cov(\hat{c}, \hat{\beta}_2)(1/\hat{\beta}_2, -\hat{c}/\hat{\beta}_2^2)^t)\}^{1/2}$. As illustrated in our simulations, the value of the standard error obtained by this approach is similar to the one obtained by the formula supplied in Proposition 1. Hence, from a practical viewpoint, the proposed method is readily implemented in software packages such as SAS and R.

## 5. Simulation studies

Simulation studies were conducted to examine the finite sample performance of the proposed estimating procedure. Specifically, normal data with identity link and binary data with logit link incorporating a change point were considered in the simulations.

For the general linear model, we assume that the response variable is normally distributed with constant variance $\sigma^2 = 1$, and that the independent variable X follows a uniform distribution with a range of $[-2, 2]$. We generated the data according to

$$\mu(x) = E(Y|x) = 2 + 3x - 5(x - 0.5)_+.$$

In this model, the change point in which the relationship between Y and x changes is set at $\tau = 0.5$. For each of the 1000 replications, we generated a sample of n=500 independent observations from the distribution $N(\mu(x), 1)$. The mean, median and standard deviation of $\hat{\delta}$ proposed in section 2 as well as their average standard errors calculated based on Proposition 1 are presented in Table 1. As discussed in Section 4, results of the estimated values as well as the precision of the estimates are not sensitive to the choice of the smoothing function. The smoothing function used in this simulation is $K(u) = 1/(2\pi)^{1/2} \int_{-\infty}^{u} e^{-t^2/2}dt$, and the smoothing parameter $h_n$ is set to be $n^{-2}$. The average s.e. by Delta method in the table is referring to the standard error obtained via the formula discussed in section 4. Note that the mean and median of the estimates are very close to the true values. This, along with the fact that the average standard errors are very similar to the corresponding sample standard deviation of the estimates suggests that our estimators converge to the true values and the normal approximation for the distribution of the estimators is valid.

The empirical coverage probabilities for $\delta$, which in this simulation are the proportion of the 95% confidence intervals (calculated by the usual estimate $\pm 1.96$se)

TABLE 1

*Estimates of the change point and the $\beta's$ from 1000 simulated samples of 500 observations for normal response variable with identity link*

| Parameters | $\beta_0$ | $\beta_1$ | $\beta_2$ | $\tau$ |
|---|---|---|---|---|
| True Value | 2 | 3 | −5 | 0.5 |
| Mean | 1.980 | 2.983 | −4.977 | 0.503 |
| Median | 1.978 | 2.983 | −4.972 | 0.505 |
| S.D. | 0.074 | 0.076 | 0.208 | 0.039 |
| Average s.e. by Proposition 1 | 0.082 | 0.078 | 0.188 | 0.037 |
| Average s.e. by Delta method | 0.076 | 0.074 | 0.202 | 0.040 |
| Empirical coverage probability(%) of normal CI | 95.9 | 95.4 | 92.8 | 94.2 |
| Empirical coverage probability of bootstrapped CI (%) | 92.6 | 93.3 | 93.0 | 94.4 |

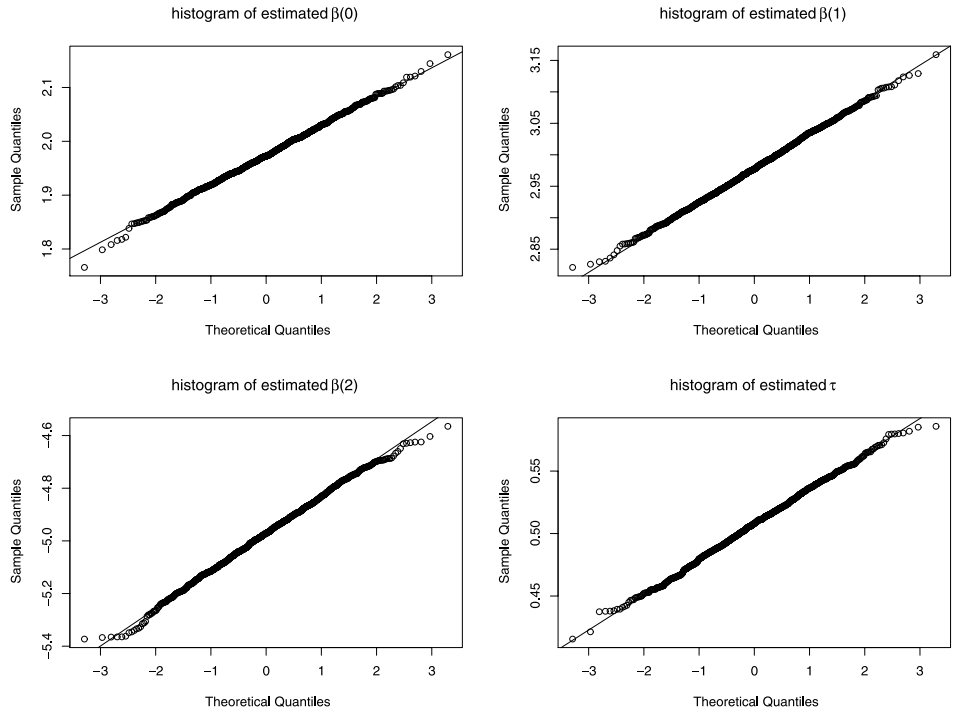

Fig 1. *Q-Q plots of the sampling distribution of the estimates from Table 1.*

for the 100 random samples containing $\delta$, by and large, are close to the nominal 95% level. In addition, the bootstrapped confidence intervals are also obtained following the bootstrap resampling scheme suggested by Boukai [4]. For each of the random samples $(X_i, Y_i), i = 1, \ldots, n$, first the change point $\tau$ as well as other parameters were estimated according to the proposed estimation procedure. The bootstrapping resample $(X_i^*, Y_i^*), i = 1, \ldots, n$ was then generated in such a way that $(X_i^*, Y_i^*), i \leq [\hat{\tau}]$ are from $(X_i, Y_i), i \leq [\hat{\tau}]$, and the rest $(X_i^*, Y_i^*), [\hat{\tau}] < i \leq n$ are from $(X_i, Y_i), [\hat{\tau}] < i \leq n$, where $[a]$ is defined as the nearest integer to a. We obtained 1000 such resamples for each of the 1000 original random samples, and then the empirical 95% bootstrap confidence interval was obtained following the application of the proposed estimation procedure to these resamples. The proportions of these 1000 bootstrap confidence intervals containing the true parameter values of $\delta$ are listed in Table 1 as the coverage probabilities of the bootstrapped CI's, which are reasonably close to the coverage probabilities of the normal CI and hence the nominal 95% level.

The Q-Q plots for the estimated $\boldsymbol{\beta}$ and the change point $\tau$ in Figure 1 also show that the sampling distributions of these estimators are close to the normal distributions. Results in Table 1 and Figure 1 suggest that the asymptotic normality of the estimators is well supported by the findings via this simulation.

For the logistic regression, we generated the regressor X from a uniform distribution $U(-2, 2)$, and the binary response variable Y according to

$$pr(Y = 1|x) = 1/[1 + exp\{-(2 + 3x - 5(x - 0.5)_+)\}].$$

As in the normal response case, we also generated a random sample of 500 from this model with 1000 replications. The smoothing parameter used here is $h_n = n^{-3}$.

*Estimates of the change point and the $\beta'$s from 1000 simulated samples of 500 observations for the binary response variable with logit link*

| Parameters | $\beta_0$ | $\beta_1$ | $\beta_2$ | $\tau$ |
|---|---|---|---|---|
| True Value | 2 | 3 | −5 | 0.5 |
| Mean | 2.029 | 3.036 | −5.058 | 0.490 |
| Median | 2.014 | 3.023 | −5.037 | 0.487 |
| S.D. | 0.282 | 0.333 | 0.647 | 0.146 |
| Average s.e. by Proposition 1 | 0.273 | 0.322 | 0.668 | 0.145 |
| Average s.e. by Delta method | 0.280 | 0.329 | 0.642 | 0.154 |
| Empirical coverage probability(%) of normal CI | 94.9 | 95.5 | 95.3 | 96.0 |
| Empirical coverage probability of bootstrapped CI (%) | 93.0 | 94.5 | 93.8 | 94.9 |

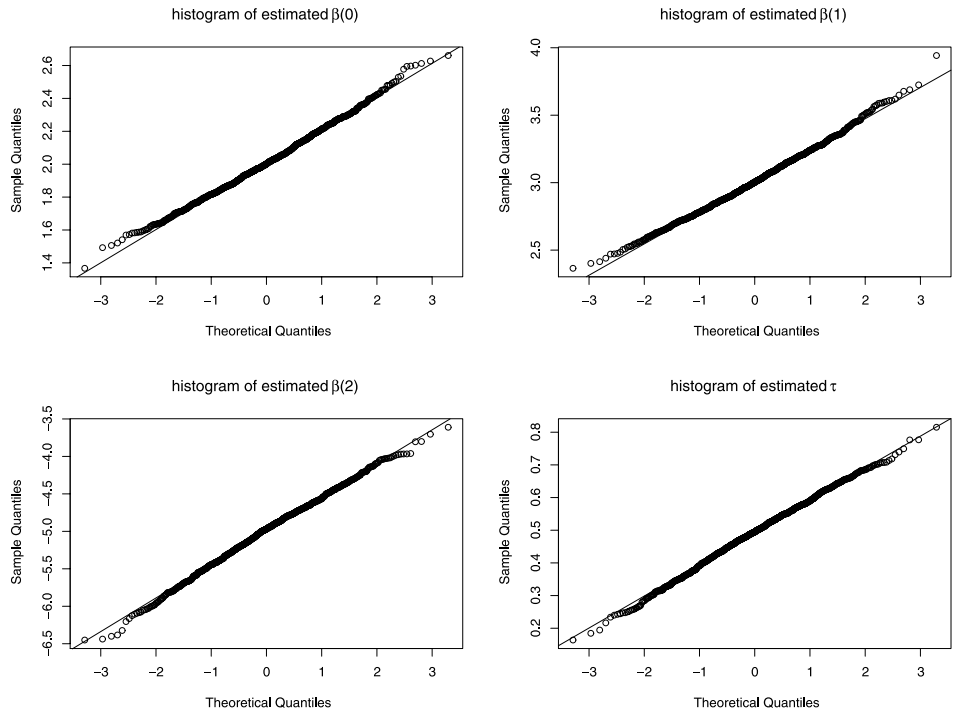

Fɪɢ 2. *Q-Q plots of the sampling distribution of the estimates from Table 2.*

Results as shown in Table 2 and Figure 2 suggest that similar conclusions can be obtained for binary data as well.

## 6. An example: The EURAMIC study

In this section, we applied, as an illustration, the proposed estimating method in Section 2 to the data from the EURAMIC study (EURopean study on Antioxidants, Myocardial Infarction, and Breast Cancer). The EURAMIC study (Kardinaal et al. [12]) is an international case-control study conducted in eight European countries and Israel, which was designed primarily to evaluate the association of antioxidants with the risk of developing a first myocardial infarction (MI) in men aged older than 70. Our example focuses on the portion of the data involving the dose-response

TABLE 3

*Parameter estimates and their standard error estimates (when available) and 95% confidence intervals relating alcohol intake to the risk of MI in the EURAMIC study*

| Parameter (s.e.) (95% C.I.) | $\beta_0$ | $\beta_1$ | $\beta_2$ | $\tau$ |
|---|---|---|---|---|
| Pastor and Guallar [16] | | 0.008 | 0.009 | 13.118 |
| | | (0.000, 0.021) | (0.001, 0.084) | (4.679, 50.062) |
| | −11.64 | 0.008 | 0.009 | 13.10 |
| Proposed method | (1.49) | (0.004) | (0.007) | (4.68) |
| | (−14.57, −8.71) | (0.0002, 0.016) | (−0.005, 0.022) | (3.94, 22.27) |

of alcohol intake and risk of myocardial infarction. For comparison purposes, we used the same sample as in Pastor and Guallar [16], who entertained the change point model. In this sample, there are 330 cases who had a confirmed diagnosis of first acute MI and 441 controls who were obtained through several random sample schemes. The primary risk factor is alcohol intake during the year before the study took place. It is well-recognized in the literature that the risk of MI as a function of alcohol intake is J-shaped. To capture this phenomenon, Pastor and Guallar [16] have fitted a quadratic-linear model to relate the risk of MI, in logit scale, to alcohol intake, adjusting for other covariates including age, smoking status, waist-hip ratio, history of diabetes, history of hypertension, family history of coronary disease and the dummy variables identifying the medical centers where the data were obtained.

Table 3 shows the estimates and their confidence intervals reported by Pastor and Guallar [16]. The approach they took is likelihood ratio-based using the original likelihood which has only the continuous first derivative with regard to $\tau$. In their approach, the authors estimated the regression coefficients and $\tau$ by maximizing the likelihood function without smoothing. Recognizing the breakdown of the conventional asymptotic results, they computed the likelihood ratio-based confidence intervals for $\tau$ and other parameters.

Also presented in Table 3 are results using the proposed estimating method. Here $K(\cdot)$ is chosen to be the cumulative distribution function of the standard normal random variable, and $h_n = n^{-3}$, where $n = 771$, the total sample size for this data set. Both approaches give rise to very similar point estimates for the change point $\tau$ as well as the $\boldsymbol{\beta}$ coefficients. Specifically, the threshold for which the risk of MI changes its direction is estimated at $13.1 (\pm 4.68)$ grams per day. However, the confidence intervals of $\tau$ are rather different from each other. While the approach by Pastor and Guallar [16] did not completely address the problem resulting from the nonsmoothness of the likelihood function with respect to the change point parameter, our estimated confidence interval for $\tau$ ($\hat{\tau} \pm 1.96 s.e.$) is based on a valid asymptotic theory and, at least for this example, provides a much tighter interval for investigators to pinpoint the point where the well-known protective alcohol intake effect may be reversed.

## 7. Discussion

Many approaches have been suggested to capture the phenomenon of abrupt change in relating the response variable to a particular independent variable. These include non-parametric smoothing and transforming or categorizing the continuous independent variables, etc. Like other approaches, the change point model is unlikely to fully capture the underlying mechanism. This approach, however, is appealing in that a specific parameter, $\tau$ in this case, is introduced to quantify the scientific objective of interest. In the alcohol intake and risk of MI example discussed above, by estimating the threshold for alcohol intake where the risk of MI is apparently heightened and the direction altered, the clinician is in a position to inform patients at risk the maximum allowable level of daily alcohol consumption. In fact, clinical advice as such is given daily to patients at risk for various diseases. Therefore, obtaining more precise estimates for these threshold values via the change point models rather than giving a simple cutoff point via conventional wisdom or experience would have significant clinical implications for health care practice.

On the other hand, the introduction of the change point parameter into the statistical model brings intriguing theoretical difficulties in both detection and estimation of such phenomenon. Our primary objective in this paper is to provide an estimate for the change point with desirable asymptotic properties under the GLM framework. The consistency and the asymptotic normality of the proposed estimators, whose variance can be easily estimated, enhances the opportunity for researchers to make statistical inference on the change point parameter.

Our approach of using the modified objective function eliminates the nonsmoothness problem with the change point parameter in the likelihood function. Similar modification can be applied to the situation where only the first two moments of the response variable is available, and the full knowledge of the probabilistic mechanism is absent. This work will be reported elsewhere.

## Appendix A: Proof of consistency

First, we state the consistency of the estimators, which is needed in the proof asymptotic normality, as

**Proposition A.1** (Consistency). *Let $\hat{\boldsymbol{\delta}} = (\hat{\boldsymbol{\beta}}, \hat{\tau})$ be the estimators which maximize the objective function $Q_n(\boldsymbol{\beta}, \tau)$ and $\boldsymbol{\delta^0}$ be the true value of $\boldsymbol{\delta}$, then under some regularity conditions and given (i),(ii) and (iii), $\hat{\boldsymbol{\delta}}$ converges in probability to $\boldsymbol{\delta}^0$, i.e., $\hat{\boldsymbol{\delta}} \xrightarrow{P} \boldsymbol{\delta}^0 = (\boldsymbol{\beta}^0, \tau^0)$.*

To show the consistency of the estimators maximizing the objective function $Q_n(\beta, \tau)$, we need the following lemma.

**Lemma A.1.** *For a sequence of random variables $Y_i$ satisfying (2.2), with the same conditions as in Proposition A.1,*

(A) $\sup\limits_{(\boldsymbol{\beta}, \tau)} \dfrac{1}{n} |Q_n(\boldsymbol{\beta}, \tau) - l_n(\boldsymbol{\beta}, \tau)| \xrightarrow{P} 0.$

*Proof.* First, under the condition $\lim\limits_{n \to \infty} P(|X_n - \tau| < \epsilon) = 0, \epsilon > 0$, there exists a $N > 0$, s.t. for $n > N, P(|X_n - \tau| < \epsilon) < \eta$, for any $\eta > 0$. $K(\frac{x_i - \tau}{h_n}) \to I(x_i > \tau)$ uniformly for $i > N$. In addition, if $|X_i - \tau| \geq \epsilon$ and $K(\cdot)$ satisfies equation (4.1), it can be easily shown

(A.1) $$K(\frac{x_i - \tau}{h_n}) = I(x_i > \tau) + o(h_n^m).$$

Note that $b'(\cdot) = g^{-1}(\cdot)$,

$$
\begin{aligned}
D &= \frac{1}{n}(Q_n(\boldsymbol{\beta}, \tau) - l_n(\boldsymbol{\beta}, \tau)) \\
&= \frac{1}{n} \sum_{i=1}^n (Y_i - g^{-1}(\theta^*)) \beta_2 (x_i - \tau)(K(\frac{x_i - \tau}{h_n}) - I(x_i > \tau)),
\end{aligned}
$$

where $|\theta_i^* - (\beta_0 + \beta_1 + \beta_2(x_i - \tau)_+)| < |\beta_2(K(\frac{x_i - \tau}{h_n}) - I(x_i > \tau))|$.

$$
\begin{aligned}
D \;=\; & \frac{1}{n}\sum_{i=1}^{N}(Y_i - g^{-1}(\theta^*))\beta_2(x_i - \tau)(K(\frac{x_i - \tau}{h_n}) - I(x_i > \tau))I(|x_i - \tau| < \epsilon) \\
+ \;& \frac{1}{n}\sum_{i=N+1}^{n}(Y_i - g^{-1}(\theta^*))\beta_2(x_i - \tau)(K(\frac{x_i - \tau}{h_n}) - I(x_i > \tau))I(|x_i - \tau| < \epsilon) \\
+ \;& \frac{1}{n}\sum_{i=1}^{n}(Y_i - g^{-1}(\theta^*))\beta_2(x_i - \tau)(K(\frac{x_i - \tau}{h_n}) - I(x_i > \tau))I(|x_i - \tau| \geq \epsilon).
\end{aligned}
$$

Since $var(Y_i) < \infty$, by the weak law of large numbers,

$$
\frac{1}{n}\sum_{i=1}^{n}[Y_i - E(Y_i)] \xrightarrow{P} \lim_{n\to\infty}\frac{1}{n}\sum_{i=1}^{n}E[Y_i - E(Y_i)] < \infty.
$$

With some constant $C > 0$, combining this with the condition $P(|X_n - \tau| < \epsilon) \to 0$ as $n \to \infty$,

$$
\begin{aligned}
\sup_{(\boldsymbol{\beta},\tau)}\frac{1}{n}|Q_n(\boldsymbol{\beta},\tau) - l_n(\boldsymbol{\beta},\tau)| \;\leq\; & \sup_{(\boldsymbol{\beta},\tau)}\frac{1}{n}|\sum_{i=1}^{N}[Y_i - g^{-1}(\theta^*)]| \cdot 2C \cdot 1 \\
& + \sup_{(\boldsymbol{\beta},\tau)}\frac{1}{n}|\sum_{i=N+1}^{n}[Y_i - g^{-1}(\theta^*)]| \cdot 2C \cdot \eta \\
& + \sup_{(\boldsymbol{\beta},\tau)}\frac{1}{n}|\sum_{i=1}^{n}[Y_i - g^{-1}(\theta^*)]| \cdot C \cdot o(h_n^m) \cdot 1.
\end{aligned}
$$

Since N is finite, the first term converges to 0 as $n \to \infty$, and $\eta$ can be made arbitrarily small, hence we have

$$
\sup_{(\boldsymbol{\beta},\tau)}\frac{1}{n}|Q_n(\boldsymbol{\beta},\tau) - l_n(\boldsymbol{\beta},\tau)| \xrightarrow{P} 0. \qquad\qquad \square
$$

*Proof of Proposition A.1.* The log likelihood function for the observed $Y's$ is

$$
\begin{aligned}
l_n(\boldsymbol{\beta},\tau) \;=\; & \sum_{i=1}^{n}log f(Y_i; \boldsymbol{\beta},\tau,x_i) \\
=\; & \sum_{i=1}^{n}[Y_i\{\beta_0 + \beta_1 x_i + \beta_2(x_i - \tau)_+\} \\
& - b(\beta_0 + \beta_1 x_i + \beta_2(x_i - \tau)_+) + log c(Y_i)].
\end{aligned}
$$

By the weak law of large numbers, $\frac{1}{n}\sum_{i=1}^{n}[log f(Y_i; \boldsymbol{\beta},\tau,x_i) - E\{log f(Y_i; \boldsymbol{\beta},\tau,x_i); \beta^0, \tau^0\}] \xrightarrow{P} 0$.

Let

$$
l_0(\boldsymbol{\beta},\tau) = \lim_{n\to\infty}\frac{1}{n}\sum_{i=1}^{n}E\{log f(Y_i; \boldsymbol{\beta},\tau,x_i); \beta^0, \tau^0\}.
$$

By Jensen's inequality, and since

$$E\{\frac{exp(Y(\beta_0 + \beta_1 x + \beta_2(x-\tau)_+) - b(\beta_0 + \beta_1 x + \beta_2(x-\tau)_+) + c(Y))}{exp(Y(\beta_0^0 + \beta_1^0 x + \beta_2^0(x-\tau^0)_+) - b(\beta_0^0 + \beta_1^0 x + \beta_2^0(x-\tau^0)_+) + c(Y))};$$
$$\beta^0, \tau^0\} = 1,$$

then

$$E\{-log(\frac{exp(Y(\beta_0 + \beta_1 x + \beta_2(x-\tau)_+) - b(\beta_0 + \beta_1 x + \beta_2(x-\tau)_+) + c(Y))}{exp(Y(\beta_0^0 + \beta_1^0 x + \beta_2^0(x-\tau^0)_+) - b(\beta_0^0 + \beta_1^0 x + \beta_2^0(x-\tau^0)_+) + c(Y))});$$
$$\beta^0, \tau^0\}$$
$$> 0.$$

Hence we have:

(B) $l_0(\boldsymbol{\beta}, \tau)$ maximizes at the true parameter values $(\boldsymbol{\beta}^0, \tau^0)$.

(C) If $X's$ are bounded, and $(\boldsymbol{\beta}, \tau) \in$ is a compact set, it follows that since $Y(\beta_0 + \beta_1 X + \beta_2(x-\tau)_+) - b(\beta_0 + \beta_1 x + \beta_2(x-\tau)_+) \leq M|Y| + N$ and $E[|Y|] < \infty$, similarly, as in the proof for consistency of maximum likelihood estimator, we have

(D) $\sup_{(\boldsymbol{\beta}, \tau)} |\frac{1}{n} l_n(\boldsymbol{\beta}, \tau) - l_0(\boldsymbol{\beta}, \tau)| \xrightarrow{P} 0$.

(E) $l_0(\boldsymbol{\beta}, \tau)$ is continuous in $(\boldsymbol{\beta}, \tau)$.

By (A), (D) and the triangle inequality,

$$\sup_{(\boldsymbol{\beta}, \tau)} |\frac{1}{n} Q_n(\boldsymbol{\beta}, \tau) - l_0(\boldsymbol{\beta}, \tau)| \quad \leq \quad \sup_{(\boldsymbol{\beta}, \tau)} |\frac{1}{n}(Q_n(\boldsymbol{\beta}, \tau) - l_n(\boldsymbol{\beta}, \tau))|$$
$$+ \sup_{(\boldsymbol{\beta}, \tau)} |\frac{1}{n} l_n(\boldsymbol{\beta}, \tau) - l_0(\boldsymbol{\beta}, \tau)| \xrightarrow{P} 0$$

This implies that $Q_n(\beta, \tau)$ satisfies the hypothesis of Theorem 4.1.1 of Amemiya (1985). Hence we have $\hat{\boldsymbol{\delta}} \xrightarrow{P} \boldsymbol{\delta}^0$. □

## Appendix B: Proof of asymptotic normality

In this appendix, we prove Proposition 1. We first introduce the expressions of Hessian matrix elements and some lemmas needed to prove Proposition 1. Aside from the conditions stated in the lemmas, the same regularity conditions as in Proposition 1 are implied as well.

For notational purposes, let $\theta_i = \beta_0 + \beta_1 x_i + \beta_2(x_i - \tau)K(\frac{x_i - \tau}{h_n})$, then

$$S_n(\delta) = \sum_{i=1}^n \frac{\partial q_i(\delta)}{\partial \theta_i} \frac{\partial \theta_i}{\partial \delta}$$

and

$$J_n(\delta) = -\sum_{i=1}^n \{\frac{\partial^2 q_i(\delta)}{\partial \theta_i^2} \frac{\partial \theta_i}{\partial \delta} (\frac{\partial \theta_i}{\partial \delta})^t + \frac{\partial q_i(\delta)}{\partial \theta_i} \frac{\partial^2 \theta_i}{\partial \delta^2}\},$$

where

$$\frac{\partial q_i(\delta)}{\partial \theta_i} = y_i - g^{-1}(\theta_i), \qquad \frac{\partial^2 q_i(\delta)}{\partial \theta_i^2} = -(g^{-1})'(\theta_i),$$

$$\frac{\partial \theta_i}{\partial \delta} = (1, x_i, (x_i - \tau)K(\frac{x_i - \tau}{h_n}), -\beta_2(K(\frac{x_i - \tau}{h_n}) + \frac{x_i - \tau}{h_n}K'(\frac{x_i - \tau}{h_n})))^t,$$

and the nonzero elements of the matrix $\frac{\partial^2 \theta_i}{\partial \delta^2}$ are

$$\frac{\partial^2 \theta_i}{\partial \delta^2}[3,4] = \frac{\partial^2 \theta_i}{\partial \delta^2}[4,3] = -\{K(\frac{x_i - \tau}{h_n}) + \frac{x_i - \tau}{h_n}K'(\frac{x_i - \tau}{h_n})\},$$

$$\frac{\partial^2 \theta_i}{\partial \delta^2}[4,4] = \beta_2\{\frac{2}{h_n}K'(\frac{x_i - \tau}{h_n}) + \frac{x_i - \tau}{h_n^2}K''(\frac{x_i - \tau}{h_n})\}.$$

**Lemma B.1.** *Suppose that*

(a) $\sup_u K'(u) \le M < \infty$, *and* $|u|K'(u) \to 0$ *as* $|u| \to \infty$,
(b) $sup_u K''(u) \le M < \infty$, *and* $u^2 K''(u) \to 0$ *as* $|u| \to \infty$,

*and under conditions (I), (II) and (III), matrix* $\frac{1}{n}J_n$ *converges element-wise in probability to*

$$\lim_{n\to\infty}\frac{1}{n}\Sigma_n = \Sigma.$$

*Proof.* Note that there are basically three different types of terms that need to be considered in both matrices $J_n/n$ and $\Sigma_n/n$. First, $\frac{1}{n}\sum_{i=1}^{n}(g^{-1})'(\theta_i)x_i^k K(\frac{x_i-\tau}{h_n})$ for $k = 0, 1, 2$. Under the conditions (i), (ii) and condition (III) on the link function $g(\cdot)$ along with equations (2.3) and (2.4), the above terms converge to finite numbers as $n \to \infty$, with the same arguments in proof of Lemma A.1.

Similarly, by condition (a) and the similar argument in Lemma A.1, as $n \to \infty$, the term

$$\frac{1}{n}\sum_{i=1}^{n}(g^{-1})'(\theta_i)\frac{x_i-\tau}{h_n}K'(\frac{x_i-\tau}{h_n}) \to 0.$$

Lastly, note that for the natural link function, most elements of matrices $\Sigma_n$ and $J_n$ are identical and free of $Y's$ except elements $J_n[3,4], J_n[4,3]$ and $J_n[4,4]$, which are $\partial^2 Q_n(\boldsymbol{\beta},\tau)/\partial\beta_2\partial\tau, \partial^2 Q_n(\boldsymbol{\beta},\tau)/\partial\tau\partial\beta_2$ and $\partial^2 Q_n(\boldsymbol{\beta},\tau)/\partial\tau^2$, respectively. Both the terms in $\partial^2 Q_n(\boldsymbol{\beta},\tau)/\partial\beta_2\partial\tau$ and in $\partial^2 Q_n(\boldsymbol{\beta},\tau)/\partial\tau^2$ involving $Y's$ can be shown to approach 0 as $n \to \infty$. Therefore, $\frac{1}{n}\Sigma_n$ and $\frac{1}{n}J_n$ converge to the same matrix $\Sigma$. $\qquad\square$

**Lemma B.2.** *Under conditions (I) and (II), the normed estimating function is asymptotically normal,*

$$\Sigma_n^{-1/2}\boldsymbol{S_n} \overset{D}{\sim} N(\boldsymbol{0}, I_4).$$

*Proof.* For fixed $\omega > 0$ and the unit vector $\boldsymbol{\lambda}, \boldsymbol{\lambda}'\boldsymbol{\lambda} = 1$, we have the sequence $\boldsymbol{\delta_n} = \boldsymbol{\delta^0} + \omega(\Sigma_n^{-1/2})^t\boldsymbol{\lambda} \in N_n(\omega)$. The Taylor expansion of the objective function is

$$Q_n(\boldsymbol{\delta_n}) = Q_n(\boldsymbol{\delta^0}) + (\boldsymbol{\delta_n} - \boldsymbol{\delta^0})'\boldsymbol{S_n} - (\boldsymbol{\delta_n} - \boldsymbol{\delta^0})'J_n(\bar{\boldsymbol{\delta}}_n)(\boldsymbol{\delta_n} - \boldsymbol{\delta^0})/2,$$
$$||\bar{\boldsymbol{\delta}} - \boldsymbol{\delta^0}|| < ||\boldsymbol{\delta}_n - \boldsymbol{\delta^0}||.$$

Taking exponential and rearranging,

$$exp\{(\boldsymbol{\lambda}'V_n(\bar{\boldsymbol{\delta}}_n)\boldsymbol{\lambda}\omega^2/2) + Q_n(\boldsymbol{\delta_n})\} = exp\{(\omega\boldsymbol{\lambda}'\Sigma_n^{-1/2}\boldsymbol{S_n}) + Q_n(\boldsymbol{\delta^0})\}.$$

Hence,

(B.1)
$$exp\{(\boldsymbol{\lambda}'V_n(\bar{\boldsymbol{\delta}}_n)\boldsymbol{\lambda}\omega^2/2)\}\frac{exp\{Q_n(\boldsymbol{\delta_n})\}}{exp\{l_n(\delta_n)\}}L_n(\delta_n)$$
$$= exp\{(\omega\boldsymbol{\lambda}'\Sigma_n^{-1/2}\boldsymbol{S_n})\}\frac{exp\{Q_n(\boldsymbol{\delta^0})\}}{exp\{l_n(\delta^0)\}}L_n(\delta^0),$$

where $L_n(\cdot)$ denotes the likelihood function of Y. By Lemma A.1, $\sup\limits_{(\boldsymbol{\beta},\tau)} \frac{1}{n}|Q_n(\boldsymbol{\delta}) - l_n(\boldsymbol{\delta})| \to 0$, and by condition (II), we have for any $\epsilon > 0$, there exists $n_1 > 0$, such that for $n \geq n_1$,

$$|exp(\boldsymbol{\lambda}'V_n(\tilde{\boldsymbol{\delta}_n})\boldsymbol{\lambda}\omega^2/2) - exp(\omega^2/2)| < \epsilon,$$

hence

$$|E[exp(\boldsymbol{\lambda}'V_n(\tilde{\boldsymbol{\delta}_n})\boldsymbol{\lambda}\omega^2/2)] - exp(\omega^2/2)|$$
$$\leq E[|exp(\boldsymbol{\lambda}'V_n(\tilde{\boldsymbol{\delta}_n})\boldsymbol{\lambda}\omega^2/2) - exp(\omega^2/2)|] < \epsilon.$$

The expectation of the left side of (B.1) converges to the moment-generating function of the standard normal. It follows that the expectation of the right side of (B.1) also converges to $exp(\omega^2/2)$. Thus, $\boldsymbol{\lambda}'\Sigma_n^{-1/2}\boldsymbol{S_n}$ is asymptotically standard normal, that is

$$\Sigma_n^{-1/2}S_n \overset{D}{\sim} N(\boldsymbol{0}, I_4).$$

Finally, proof of Proposition 1 completes with

$$0 = \boldsymbol{S_n}(\hat{\boldsymbol{\delta}}) = \boldsymbol{S_n} - J_n(\boldsymbol{\delta}^*)(\hat{\boldsymbol{\delta}} - \boldsymbol{\delta}^0), \quad ||\boldsymbol{\delta}^* - \boldsymbol{\delta}^0|| < ||\hat{\boldsymbol{\delta}} - \boldsymbol{\delta}^0||,$$

Lemma B.1 along with condition (II), and consistency of $\hat{\delta}$ stated in Proposition A.1. □

**Acknowledgments.** The authors wish to thank Drs Eliseo Guallar and Robert Pastor for kindly making The EURAMIC (EURopean Study on Antioxidant, Myocardial Infarction, and Breast Cancer) data set available to us.